\documentclass{amsart}
\usepackage[utf8]{inputenc}
\usepackage{enumerate}
\usepackage{enumitem}
\usepackage{graphicx}
\usepackage{wrapfig}
\usepackage{subcaption}
\usepackage{listings}
\usepackage{color}
\usepackage[foot]{amsaddr}
\definecolor{mygreen}{RGB}{28,112,30} 
\definecolor{mylilas}{RGB}{170,55,241}
\definecolor{myblue}{RGB}{20,30,171}
\definecolor{myred}{RGB}{200,10,30}
\usepackage{microtype}
\usepackage{multirow}
\usepackage{array}
\newcounter{rowcount}
\usepackage{amsmath}
\usepackage{amsfonts}
\usepackage{amssymb}
\usepackage{amscd}
\usepackage{amstext}
\usepackage{amsbsy}
\usepackage{amsopn}
\usepackage{amsthm}
\usepackage{thmtools}
\usepackage{amsxtra}
\usepackage{float}
\usepackage{bbm}
\usepackage{url}
\usepackage[numbers,square]{natbib}
\usepackage{hyperref}
\hypersetup{
colorlinks=true,
linkcolor=myred,
citecolor=blue}

\newtheoremstyle{mytheoremstyle} 
  {12pt}                    
  {\topsep}                    
  {\slshape}                   
  {}                           
  {\bfseries}                   
  {.}                          
  {.5em}                       
  {}  

\newtheoremstyle{mydefstyle} 
  {12pt}                    
  {12pt}                    
  {}                   
  {}                           
  {\bfseries}                   
  {.}                          
  {.5em}                       
  {}  

  \newtheoremstyle{plainsl}%
  	{}
  	{\topsep}
  	{\slshape} 
  	{}
  	{\normalfont\bfseries}
  	{.}
  	{ }
  	{}

\theoremstyle{plainsl}
\newtheorem{thm}{Theorem}[section]
\newtheorem*{thm*}{Theorem}

\newtheorem{cor}[thm]{Corollary}
\newtheorem{lem}[thm]{Lemma}
\newtheorem*{lem*}{Lemma}

\theoremstyle{remark}

\theoremstyle{mydefstyle}

\newtheorem*{note*}{Remark}

\newcommand{\R}{\mathbb{R}}

\newcommand{\Cx}{\mathbb{C}}
\newcommand{\A}{\mathcal{A}}
\newcommand{\B}{\mathcal{B}}
\newcommand{\C}{\mathcal{C}}

\newcommand{\E}{\mathcal{E}}

\renewcommand{\phi}{\varphi}

\newcommand{\sm}{\backslash}

\renewcommand\qed{%
	\ifmmode\eqno\sqr53
	\else\nolinebreak\ \hfill\sqr53\medbreak\fi}
\renewcommand\proof{\noindent\textsl{Proof. }}
\newcommand\sqr[2]{{\vbox{\hrule height.#2pt
    \hbox{\vrule width.#2pt height#1pt \kern#1pt
        \vrule width.#2pt}\hrule height.#2pt}}}

\def\pmat#1{{\begin{pmatrix}#1\end{pmatrix}}}

\DeclareMathOperator\sgn{sgn}

\DeclareMathOperator\Cay{Cay}

\title{A note on eigenvalues of Cayley graphs}
\author{Arnbjörg Soffía Árnadóttir $^{1}$}
\address{$^{1}$ Corresponding author. Department of Applied Mathematics and Computer Science, Technical University of Denmark, DK-2800 Lyngby, Denmark. e-mail:  {\tt sofar@dtu.dk}\\ \normalfont{Supported by Independent Research Fund Denmark, 8021-00249B AlgoGraph and Carlsberg Semper Ardens Accelerate CF21-0682 Quantum Graph Theory.}}

\author{Chris Godsil $^2$}
\address{$^2$ Department of Combinatorics and Optimization, University of Waterloo, 200 University Avenue West, Waterloo, ON, Canada N2L 3G1. e-mail:  {\tt cgodsil@uwaterloo.ca}\\ \normalfont{Supported by NSERC (Canada), Grant No.\ RGPIN-9439.}}

\begin{document}
\renewcommand{\itshape}{\slshape}


\begin{abstract}
  A graph is called \textit{integral} if all its eigenvalues are integers. A Cayley graph is called \textit{normal} if its connection set is a union of conjugacy classes. We show that a non-empty integral normal Cayley graph for a group of odd order has an odd eigenvalue.
\end{abstract}

\maketitle

\section{Introduction}

Spectral graph theory refers to the study of eigenvalues and eigenvectors of matrices arising from graphs, the adjacency matrix being the most notable one. This topic has been of interest for many years since often times, the spectral properties of a graph can tell us something about the graph structure. However, the spectrum of a graph is usually not very accessible which brings us to Cayley graphs.

Spectra of Cayley graphs can be studied using various algebraic tools, such as representation theory, association schemes and, of course, group theory. In this paper, we prove a result on the eigenvalues of Cayley graphs using some of these tools, but very little graph theory.

We call a Cayley graph \textit{normal} if it has a conjugacy-closed connection set and we call a graph \textit{integral} if its adjacency matrix has only integer eigenvalues. The main result of the paper is the following.

\begin{thm*}[Theorem \ref{thm:oddeval}]
  Let $G$ be a group of odd order and let $X$ be a non-empty integral, normal Cayley graph for $G$. Then $X$ has an odd eigenvalue.
\end{thm*}

This theorem gives some insight into the spectral graph theory of Cayley graphs, but the motivation comes from a different direction entirely, namely the study of quantum walks.
A \textit{quantum walk} on a graph with adjacency matrix $A$ is given by the matrices $U_A(t)=\exp(itA)$ where $t\in\R$. The graph is said to have \textit{perfect state transfer} from vertex $u$ to vertex $v$ at time $t$ if the $uv$-entry of $U_A(t)$ has absolute value one.

Perfect state transfer was introduced by Bose in 2003 \cite{bose2003} and is of significant interest in quantum physics. A review on this by Kay can be found in \cite{kay2010}. In short, perfect state transfer is  a useful property for a graph to have, but unfortunately quite rare.

In terms of examples, Cayley graphs play an important role. Extensive results have been proved for perfect state transfer on Cayley graphs for the elementary abelian $2$-groups, \cite{bernasconi, cheung, chan2013complex} and in 2013, Ba\v si\'c gave a complete characterization of Cayley graphs for cyclic groups admitting perfect state transfer \cite{basic2013char}. In 2022, the authors of this paper generalized Ba\v si\'c's result to abelian groups with a cyclic Sylow-$2$-subgroup \cite{arnadottir2022pst} using a more group theoretic approach.


A crucial part of the characterization in \cite{arnadottir2022pst} is a lemma that is analogous to Theorem \ref{thm:oddeval} in this paper, for abelian groups. With the long term goal of generalizing the characterization even further to non-abelian groups, this theorem will be essential. Moreover, this motivation explains the somewhat restrictive condition in the theorem that the graph be integral, since this is a necessary condition for perfect state transfer to occur in our graphs.

The main tool used in this paper will be that of association schemes which will be introduced in Section \ref{sec:schemes}. In sections \ref{sec:partitions} and \ref{sec:conj} we prove some important lemmas and relate the schemes to Cayley graphs and in Section \ref{sec:proof}, we complete the proof of the main theorem.

\section{Preliminaries}\label{sec:prelim}

Let $G$ be a group with identity $e$ and let $\C$ be a subset of $G\sm \{e\}$ such that either $\C^{-1}=\C$ or $\C^{-1}\cap \C=\emptyset$. Define a graph with vertex set $G$ and arc set $\{(g,h): hg^{-1}\in \C\}$. If $\C^{-1}=\C$, this is an undirected graph (which we call a graph from now on) and if $\C^{-1}\cap \C=\emptyset$, it is a digraph. This is the \textit{Cayley graph (or Cayley digraph)} for $G$ with respect to $\C$ and is denoted by $\Cay(G,\C)$. The set $\C$ is the \textit{connection set} of this graph.  We say that the Cayley graph is \textit{normal} if $\C$ is a union of conjugacy classes of $G$.


We define a \textit{signed Cayley graph (or digraph)} as a Cayley graph, $X=\Cay(G,\C)$ together with a weight function, $\omega:\C\to \{\pm1\}$ with the added requirement if $X$ is undirected that the fibre of both $1$ and $-1$ are inverse-closed. Denote by $\C^+$ and $\C^-$ the fibres of $1$ and $-1$, respectively. We say that a signed Cayley graph is \textit{normal} if both $\C^+$ and $\C^-$ are conjugacy-closed.

Let $X$ be a graph or digraph. The \textit{adjacency matrix} of $X$ is indexed by its vertices and defined by
\[A(X)_{uv} =
\begin{cases}
    1 & \text{if } (u,v) \text{ is an arc,}\\
    0 & \text{otherwise.}
\end{cases}\]
The \textit{eigenvalues and eigenvectors of the graph (or digraph)} are the eigenvalues and eigenvectors of its adjacency matrix. Note that a graph has a symmetric adjacency matrix, so its eigenvalues are real numbers. A graph is said to be \textit{integral} if all its eigenvalues are integers.

We see that if $X=\Cay(G,\C)$, its adjacency matrix can be expressed as
\[A(X)_{gh} =
\begin{cases}
    1 & \text{if } hg^{-1}\in \C,\\
    0 & \text{otherwise.}
\end{cases}\]
If $X$ is a signed Cayley graph, we define its \textit{signed adjacency matrix} by
\[A_{\sgn}(X) := A(\Cay(G,\C^+)) - A(\Cay(G,\C^-)).\]
The spectrum of a signed Cayley graph refers to the spectrum of this signed adjacency matrix.

\section{Association schemes}\label{sec:schemes}
The tools we use in this paper come from the theory of association schemes. In this section we give an introduction into this theory, but we refer the reader to \cite[Chapter 2]{brouwer1989} and \cite[Chapter 12]{godsil-combinatorics} for more details.

Let $J$ denote the $n\times n$ all-ones matrix. An \textit{association scheme with $d$ classes} is a set of $n\times n$ matrices, $\A = \{A_0,\dots,A_d\}$ with entries in $\{0,1\}$ such that
\begin{enumerate}
  \item $A_0=I$ and $\sum_{r=0}^d A_r = J$,
  \item $A_r^T\in \A$ for all $r$,
  \item $A_rA_s=A_sA_r$ for all $r,s$,
  \item $A_rA_s$ lies in the span of $\A$ for all $r,s$.
\end{enumerate}

The span of $\A$ is a commutative algebra, $\Cx[\A]$, called the \textit{Bose-Mesner algebra} of the association scheme and any $\{0,1\}$-matrix in this algebra is a \textit{Schur idempotent} of $\Cx[\A]$. The elements in the scheme, $A_0,\dots, A_d$, are the \textit{minimal Schur idempotents} and every Schur idempotent is a sum of some minimal Schur idempotents. If $\B$ is an association scheme such that every element of $\B$ is a Schur idempotent of $\Cx[\A]$, we say that $\B$ is a \textit{subscheme} of $\A$ (equivalently,  $\B$ is a scheme whose Bose-Mesner algebra is a subalgebra of the Bose-Mesner algebra of $\A$).


We can view the Schur idempotents of an association scheme as adjacency matrices of graphs or digraphs. A graph (or digraph) whose adjacency matrix is a Schur idempotent in the Bose-Mesner algebra of an association scheme is referred to as a \textit{graph (or digraph) in the scheme}.

The set $\A$ is a basis for $\Cx[\A]$. It can be shown that the algebra has another basis, $\E = \{E_0,\dots, E_d\}$ of matrix idempotents (we call them the \textit{minimal matrix idempotents} of the scheme) that are pairwise orthogonal and sum to the identity matrix. For $r,s=0,\dots, d$ we define scalars $p_r(s), q_r(s)$ as follows:
\[A_r = \sum_{j=0}^d p_r(j)E_j\quad\text{and}\quad E_r = \frac1n\sum_{j=0}^d q_r(j)A_j.\]

Notice that since the matrices $E_r$ are pairwise orthogonal, we have for all $r,s$
\begin{align*}
    A_rE_s &= \left(\sum_{j=0}^d p_r(j)E_j\right)E_s = p_r(s)E_s
\end{align*}
and so $p_r(s)$ is an eigenvalue of $A_r$ for all $s$ and the columns of $E_s$ are eigenvectors of $A_r$ for this eigenvalue.

Define the $d\times d$ matrices $P$ and $Q$ by
\[P=(p_r(s))_{s,r},\quad Q=(q_r(s))_{s,r}.\]
We call them the {\it eigenmatrix} and {\it dual eigenmatrix} of the scheme, respectively. We further define for $r=0,\dots, d$ integers $v_r$ and $m_r$ as the row sum of $A_r$ and rank of $E_r$, respectively. We call $v_0,\dots, v_d$ the {\it valencies} of the association scheme and $m_0,\dots, m_d$ its {\it multiplicities}. Let $D_v$ and $D_m$ be the diagonal matrices with the valencies and multiplicities respectively on the diagonal. The following identities are well known (see for example \cite[Section 2.2]{brouwer1989})
\begin{align*}
    PQ &= nI, \\ D_mP &= Q^*D_v
\end{align*}
where $Q^*$ denotes the conjugate transpose of $Q$. Combining the two identities, we get
\(P^*D_mP = nD_v\)
and the next lemma follows immediately.
\begin{lem}\label{lem:det}
  Let $\A$ be an association scheme on $n$ vertices, with $d$ classes. Let $P$ be the matrix of eigenvalues of $\A$ and let $v_0,v_1,\dots,v_d$ and $m_0,m_1,\dots,m_d$ be the valencies and multiplicities, respectively. Then
  \[\det(P^*P) = n^{d+1}\prod_{j=0}^d\frac{v_j}{m_j}.\qed\]
\end{lem}
The right hand side of the equation in Lemma \ref{lem:det} is often called the \textit{frame quotient} of the scheme. We get the following corollary.

\begin{cor}\label{cor:frame}
    The frame quotient of an association scheme is an integer.
\end{cor}
\proof
    The entries of $P$ are eigenvalues of the matrices $A_0,\dots, A_d$. These matrices have integer entries and so the entries of $P$ are algebraic integers. Therefore, $\det(P^*P)$ is an algebraic integer, but by Lemma \ref{lem:det}, it is rational and so it must be an integer.
\qed

\section{Quotients \& subschemes}\label{sec:partitions}

The aim of this section is to derive a relation between the eigenmatrices of schemes and their subschemes. 

We start by defining equitable quotient matrices.
Let $M\in \Cx^{n\times n}$ be a matrix with
\[M = \pmat{M_{11}&\cdots & M_{1k} \\ \vdots & \ddots &\vdots \\ M_{k1} & \cdots & M_{kk}},\]
where the block $M_{rs}$ is an $n_r\times n_s$ matrix. Suppose that the row sum of each block is constant and let $a_{rs}$ be the row sum of $M_{rs}$. Then the $k\times k$ matrix, $(a_{ij})_{i,j}$ is called an \textsl{equitable quotient matrix} of $M$. The following is well known (see proof for instance in \cite{you2019}).
\begin{lem}\label{lem:quotient}
  If $N$ is an equitable quotient matrix of $M$, then the characteristic polynomial of $N$ divides the characteristic polynomial of $M$.\qed
\end{lem}
We can now prove the main result of this section.
\begin{thm}
    If $\A$ is an association scheme with eigenmatrix $P$ and $\B$ is a subscheme with eigenmatrix $P'$, then $P'$ is an equitable quotient matrix of $P$. In particular, the characteristic polynomial of $P'$ divides the characteristic polynomial of $P$.
\end{thm}
\proof
Let $\A=\{A_0,\dots, A_d\}$ and $\B=\{B_0,\dots, B_k\}$. Then there is a partition of $\{0,\dots, d\}$ with cells $C_0,\dots, C_k$ (where $C_0=\{0\}$) such that
\[B_r = \sum_{j\in C_r}A_j.\]

Let $E_0,\dots, E_d$ be the minimal matrix idempotents of $\A$ and $F_0,\dots, F_k$ the minimal matrix idempotents of $\B$. Since $F_r\in \Cx[\A]$ and $E_0,\dots, E_d$ is a basis, we may write
\[F_r=\sum_{j=0}^d\theta_jE_j\]
for some $\theta_0,\dots, \theta_d$. But each $\theta_s$ will be an eigenvalue of $F_r$ and since $F_r$ is idempotent, $\theta_s\in\{0,1\}$ so each $F_r$ is a sum of a subset of the $E_j$. Further, since both the $F_r$ and the $E_r$ sum to the identity matrix, this defines another partition, of $\{0,\dots, d\}$ with cells $D_0,\dots, D_k$, where
\[F_r = \sum_{j\in D_r}E_j.\]

We now order the rows and columns of the eigenmatrix $P$ of $\A$ such that
\[P = \pmat{P_{11}&\cdots & P_{1k} \\ \vdots & \ddots &\vdots \\ P_{k1} & \cdots & P_{kk}},\]
where $P_{rs}$ consists of the eigenvalues $p_j(i)$ with $j\in C_r$ and $i\in D_s$.
We will show that each $P_{rs}$ has a constant row sum.

Since $\B$ is an association scheme, we have $B_rF_s = \lambda_r(s)F_s$ for some $\lambda_r(s)$. Recall that the minimal idempotents of a scheme are pairwise orthogonal and so for any $j\in D_s$ we have
\[F_sE_j = \left(\sum_{i\in D_s}E_i\right) E_j = E_j.\]
Therefore,
\[B_rE_j = B_rF_sE_j = \lambda_r(s)F_sE_j = \lambda_r(s)E_j\]
but we also have
\[B_rE_j = \left(\sum_{i\in C_r}A_i\right) E_j = \sum_{i\in C_r}p_i(j)E_j,\]
and notice that $\sum_{i\in C_r}p_i(j)$ is the sum of the $j$-row in $P_{rs}$. We conclude that $P_{rs}$ has constant row sum, namely $\lambda_r(s).$ We also see that this is an eigenvalue of the subscheme, $\B$.

It is now clear that the matrix, $P':=(\lambda_r(s))_{s,r}$ is an equitable quotient matrix of $P$ and is the eigenmatrix of $\B.$ The theorem then follows from Lemma \ref{lem:quotient}\qed

\section{The conjugacy class scheme}\label{sec:conj}
Let $G$ be a group of order $n$ with conjugacy classes $C_0,\dots,C_d$ (where $C_0=\{e\}$). Define the $n\times n$ matrices, $A_0,\dots,A_d$ by
\[(A_r)_{gh} =
\begin{cases}
    1 &\text{if }\, hg^{-1}\in C_r,\\
    0 &\text{otherwise.}
\end{cases}\]
Then $\A :=\{A_0,\dots, A_d\}$ is an association scheme (see for example \cite[Lemma 3.3.1]{godsil2016EKR}), we call it the \textit{conjugacy class scheme} of $G$. It is easy to see that a normal Cayley graph of $G$ is precisely a graph in its conjugacy class scheme.

\begin{lem}\label{lem:detoddint}
  If $\A$ is the conjugacy class scheme of a group of odd order, with eigenmatrix $P$, then $\det(P^*P)$ is an odd integer. 
\end{lem}
\proof
    Recall that $\det(P^*P)$ is the frame quotient of $\A,$ which is an integer by Corollary \ref{cor:frame}.

    It is clear from the definition of this scheme that the row-sum of $A_r$ is the size of the conjugacy class, $C_r$, so the valencies of the scheme are the sizes of the conjugacy classes of the group.
    Therefore, if $n$ is odd, the valencies are all odd and now the result follows from Lemma \ref{lem:det}.
\qed

We now introduce an important subscheme of the conjugacy class scheme of a group.
Define a relation on the elements of $G$ as follows. We say that $g_1$ and $g_2$ are \textit{power equivalent}, and write $g_1\approx g_2$, if $\langle g_1\rangle = \langle g_2\rangle$. This is an equivalence relation on $G$ and we refer to its equivalence classes as \textit{power classes}.

It is not too hard to verify that the power classes commute with the conjugacy classes in the sense that the closure of a power class under conjugation is equal to the closure of a conjugacy class under the power relation. This gives a new relation on the group $G$ whose classes we will call \textit{PC-classes}. Each PC-class is a union of conjugacy classes and a union of power classes of the group.

Let $D_0,\dots, D_k$ be the PC-classes of $G$ with $D_0=\{e\}$ and define $n\times n$ matrices $B_0,\dots B_k$ by
\[(B_r)_{gh} =
\begin{cases}
    1 &\text{if }\, hg^{-1}\in D_r,\\
    0 &\text{otherwise.}
\end{cases}\]
It turns out that $\B=\{B_0,\dots, B_k\}$ is an association scheme and it is clearly a subscheme of the conjugacy class scheme of $G$. Furthermore, the matrices in $\Cx[\B]$ that have integer entries also have integer eigenvalues. In particular, the graphs in this scheme are integral. We call $\B$ the \textit{integral conjugacy class scheme of $G$}.

The next theorem shows that something stronger holds, namely that all the integral normal Cayley graphs live in this scheme. This was first shown for abelian groups by Bridges and Mena in 1982 \cite[Theorem 2.4]{Bridges1982}. In 2012, Alperin and Peterson proved one direction for groups in general \cite[Theorem 4.1]{Alperin2012}, and finally the theorem was proved in 2014 by Godsil and Spiga in 2014\cite{godsil2014rationality}.

\begin{thm}{\cite[Theorem 1.1]{godsil2014rationality}}
    A graph is an integral, normal Cayley graph for a group $G$ if and only if it lies in the integral conjugacy class scheme of $G$.\qed
\end{thm}

It is now easy to prove that the eigenmatrix of the integral conjugacy class scheme for a group of odd order has odd determinant.

\begin{lem}
  Let $G$ be a group of odd order, and let $\A$ be the integral conjugacy class scheme of $G$ with eigenmatrix $P$. Then $\det(P)$ is an odd integer.
\end{lem}
\proof
  Firstly, it is clear that $\det(P)$ is an integer since the entries of $P$ are integers. Therefore, we have $\det(P^*P)=\det(P)^2$. Let $P'$ be the eigenmatrix of the conjugacy class scheme.
  Then $\det(P)$ divides $\det(P')$ and so $\det(P)^2 = \det(P^*P)$ divides $\det(P'^*P')$ which by Lemma \ref{lem:detoddint} is an odd integer. It follows that $\det(P)^2$ is an odd integer, and therefore so is $\det(P)$.
\qed

\section{Main Theorem}\label{sec:proof}

We are ready to state and prove our main theorem.

\begin{thm}\label{thm:oddeval}
  Let $G$ be a group of odd order and let $X$ be a non-empty integral, normal Cayley graph for $G$. Then $X$ has an odd eigenvalue.
\end{thm}
\proof
    Since $X$ is an integral, normal Cayley graph for $G$, it lies in the integral the conjugacy class scheme,  $\A=\{A_0,\dots, A_d\}$ of $G$. Let $P$ be the matrix of eigenvalues of this scheme. The adjacency matrix of $X$ can be written as
    \[A = \sum_{r\in C}A_r\]
    where $C\subseteq \{0,\dots, d\}$. Then, if $x$ is the characteristic vector of the subset $C$, then $Px$ is a $d+1$ vector whose entries are the eigenvalues of $X$.

    If the entries of $Px$ are all even, then $P x\equiv 0\pmod 2$, but since $\det(P)$ is odd, $P$ is invertible modulo $2$, so this implies that $X$ is empty.
\qed

The above theorem holds slightly more generally, that is, in the case where our Cayley graph is signed.

\begin{cor}
  Let $G$ be a group of odd order and let $X=\Cay(G,\C)$ be a non-empty normal signed Cayley graph for $G$ such that both $\Cay(G,\C^+)$ and $\Cay(G,\C^-)$ are integral. Then $X$ has an odd eigenvalue.
\end{cor}
\proof
We can modify the proof of Theorem \ref{thm:oddeval}, replacing $x$ by a $\{0,\pm 1\}$-vector, in the obvious way to get the eigenvalues of the signed adjacency matrix as the entries of $Px$. The rest of the proof is the same. \qed

\section{Conclusion}
We have shown that every non-empty, normal Cayley graph for a group of odd order, that has only integer eigenvalues must have an odd eigenvalue. We further showed that this holds slightly more generally, where the graph is allowed to be signed.

\section*{Acknowledgements}
The authors would like to thank David Roberson for his suggestions.


\end{document}